\documentclass[a4paper]{amsart}
\usepackage{graphicx, amsthm, amsmath, amssymb, mathrsfs, amscd}
\usepackage{amsmath}
\usepackage{breakurl}
\usepackage[tableposition=below]{caption}
\usepackage{enumitem}
\usepackage{hyperref}
\usepackage{color}
\usepackage{color, colortbl}
\usepackage{tikz-cd}

\newtheorem{theorem}{Theorem}[section]
\newtheorem{lemma}[theorem]{Lemma}
\newtheorem{cor}[theorem]{Corollary}
\newtheorem{prop}[theorem]{Proposition}

\theoremstyle{definition}
\newtheorem{rem}[theorem]{Remark}
\theoremstyle{definition}
\newtheorem{definition}[theorem]{Definition}

\newcommand{\N}{\mathbb{N}}
\newcommand{\Z}{\mathbb{Z}}
\newcommand{\Q}{\mathbb{Q}}
\newcommand{\F}{\mathbb{F}}

\title{A newform theory for Katz modular forms}

\author[Daniel Berhanu]{Daniel Berhanu}
\address{Unit$\acute{\text{e}}$ de Recherche en Math$\acute{\text{e}}$matiques\\
Universit$\acute{\text{e}}$ du Luxembourg\\
Maison du nombre\\
6, avenue de la Fonte\\
L-4364 Esch-sur-Alzette\\
Luxembourg}
\email{daniel.mamo@uni.lu\\
danielberhanu86@gmail.com}


\begin{document}

\begin{abstract}
In this paper, a strong multiplicity one theorem for Katz modular forms is studied. We show that a cuspidal Katz eigenform which admits an irreducible Galois representation is in the level and weight old space of a uniquely associated Katz newform. We also set up multiplicity one results for Katz eigenforms which have reducible Galois representation.
\end{abstract}

\maketitle

\section{Introduction}\label{Sec 1}
\setcounter{equation}{0}

The study of relations between the coefficients of classical modular forms by L. Atkin and J. Lehner in \cite{AL} led to the invention of the theory of newforms. Atkin and Lehner used the $L$-functions associated to the newforms for their investigation. W. Li in \cite{Li}, using the notion of trace operators, obtained the generalization of the Atkin-Lehner theory to the case of modular forms over congruence subgroups parameterized by two variables with characters. In this note we will generalize some of these results to Katz modular forms over $\overline{\mathbb{F}}_p$.

Katz modular forms are modular forms defined via algebraic geometry methods by N. Katz in \cite{NK}. They are defined over any ring in which the level is invertible. We work with Katz modular forms over $\overline{\mathbb{F}}_p$. Thus we always assume that the prime $p$ does not divide the levels of our modular forms. Katz modular forms admit Hecke operators analogously to holomorphic modular forms. We denote the space of Katz modular forms in level $\Gamma_1(N)$, weight $k$ and Dirichlet character $\varepsilon$ by $M_{k}(\Gamma_1(N), \varepsilon, \overline{\mathbb{F}}_p)_{\text{Katz}}$ and its cuspidal subspace by $S_{k}(\Gamma_1(N), \varepsilon, \overline{\mathbb{F}}_p)_{\text{Katz}}$. When we do not write the Dirichlet character $\varepsilon$ we assume that we did not fix any character.

Let $f$ be a normalised Katz eigenform of level $\Gamma_1(N)$ with $p\nmid N$, weight $k$ and character $\varepsilon$, with coefficients in $\overline{\mathbb{F}}_p$.  Let $f(T_l)$ be the $l$th eigenvalue of $f$ for the Hecke operator $T_l$. Then thanks to the works in \cite{DS}, \cite{Bj} there exists a unique 2-dimensional semi-simple continuous representation $\rho_f: \text{G}_{\Q} \to \text{GL}_2(\overline{\mathbb{F}}_p)$, which is unramified outside $pN$ and has the property that $\text{tr}(\rho_f(\text{Frob}_l))=f(T_l)$ and $\det (\rho_f(\text{Frob}_l))=\varepsilon(l) l^{k-1}$ for all primes $l\nmid pN$. We prefer to state the results in terms of Galois representations because they shorten the statements. However it would be possible to avoid that language for most statements.

Let us informally introduce the notation of level and weight old spaces. They are defined more precisely in the next section. First, like the classical case, we have level degeneracy maps on Katz modular forms. Let $f \in M_{k}(\Gamma_1(N), \overline{\F}_p)_{\text{Katz}}$ be any Katz modular form and let $d\geq 1$ be an integer coprime to $p$. Then we have the $d$-th degeneracy map $f(q) \mapsto f(q^d)$, which increases the level by a multiple of $d$. Then the $level~old~space$ of $f$ in the level $M$ divisible by $N$ is the $\overline{\mathbb{F}}_p$ vector space generated by modular forms $f(q^d)$ where $d$ runs through all possible divisors of $M/N$. Second, we have the following weight degeneracy maps on Katz modular forms. The map $\alpha_A$ defined by $f \mapsto Af$, where $A\in M_{p-1}(\Gamma_1(1), \overline{\F}_p)_{\text{Katz}}$ is the Hasse invariant with $q$-expansion at the cusp infinity equal to 1. It adds $p-1$ to the weight but does not change the $q$-expansion. The Frobenius map takes a form $f(q)$ to its Frobenius $\text{Frob}(f)(q)=f(q^p)$. It multiplies the weight by $p$ but does not change the level. Thus, the missing degeneracy map $q\mapsto q^p$ in the level is provided by the Frobenius.

Then by the $level~and~weight~old~space$ of $f$ in level $M$, a multiple of $N$, and weight $k'\geq k$ we understand the space generated by the images of $f$ under all possible combinations of the level and weight degeneracy maps targeted to the space of modular forms of level $M$ and weight $k'$.

We will use the following definitions of minimal levels and weights to introduce our newforms. Let $d\geq 2$ be a positive integer and let $f \in M_k(\Gamma_1(N), \overline{\mathbb{F}}_p)_{\text{Katz}}$ be any Katz eigenform. Then $f$ is said to have a $d$-$minimal~level$ if $\rho_f$ does not arise from any non-zero Katz eigenform of level $MN/d^{m}$ where $M$ and $m$ are any positive integers such that $\text{gcd}(M, d)=1$. A Katz modular form $f$ is said to have a $minimal~weight~k$ if the associated mod $p$ Galois representation $\rho_f$ does not arise from any non-zero Katz eigenform of weight strictly smaller than $k$. We assume that any Katz modular form in non-integral level is equal to zero.

\begin{definition}\label{df01}
A normalised Katz eigenform $f\in M_{k}(\Gamma_1(N), \varepsilon, \overline{\mathbb{F}}_p)_{\text{Katz}}$ is called a \emph{Katz~newform} if $N$ is an $l$-minimal level for any prime $l|N$ and $k$ is the minimal weight.
\end{definition}

The motivation behind the definition of our Katz newform is that it satisfies some of the analog results of classical newform theory.

The aim of the paper is to prove the following strong multiplicity one theorems.
\begin{theorem}\label{th02}
Let $f\in S_k(\Gamma_1(N), \varepsilon, \overline{\mathbb{F}}_p)_{\text{Katz}}$ and $g\in S_{k'}(\Gamma_1(N'), \varepsilon', \overline{\mathbb{F}}_p)_{\text{Katz}}$ be Katz newforms with $a_l(f)=a_l(g)$ for each $l$ in a set of primes of density 1. Then $f=g, k=k', N=N'$, and $\varepsilon = \varepsilon '$.
\end{theorem}

This means that Katz newforms are uniquely characterised by their associated mod $p$ Galois representations.

\begin{theorem}\label{th01} Let $F\in S_{k'}(\Gamma_1(M), \overline{\mathbb{F}}_p)_{\text{Katz}}$ be a normalised Katz eigenform. Suppose that there is a Katz newform $f\in S_{k}(\Gamma_1(N), \overline{\mathbb{F}}_p)_{\text{Katz}}$ such that $\rho_F \cong \rho_f$. Then $F$ is in the level and weight old space of $f$.
\end{theorem}

As a consequence, one can determine all possible coefficients of any normalised Hecke eigenform from its associated Katz newform, provided that it exists. See Corollary~\ref{le26} for explicit expressions.

The existence of Katz newforms is established in Theorem~\ref{pr132} (Serre's conjecture) when the associated mod $p$ Galois representation is irreducible and in Proposition~\ref{pr31} when the associated mod $p$ Galois representation is reducible of a certain type.

In Section~\ref{Sec 4}, we set up the theory of newforms for the space of Katz Eisenstein series. In the case where cuspidal Katz eigenforms have reducible mod $p$ Galois representations Eisenstein series come into the picture to describe their associated newforms. We have shown in Theorem~\ref{th47} that, under some condition, up to a suitable power multiple of the Hasse invariant, any non-ordinary cuspidal Katz eigenform with a reducible mod $p$ Galois representation is in the level old space of an associated Katz eigenform which has an optimal level.

\section{Notation and Preliminaries}\label{Sec 2}
\setcounter{equation}{0}

To state some of the theoretical results that we need later let us set the following notation. For any continuous Galois representation $\rho_p: \text{Gal}(\overline{\Q}_p/\Q_p) \to \text{GL}_2(\overline{\F}_p)$ we associate an integer $k(\rho_p)$ as in the Definition 4.3 of \cite{BE}. Then for continuous mod $p$ representation $\rho:\text{G}_{\Q}\to \text{GL}_2(\overline{\F}_p)$ we set $k(\rho)=k(\rho|_{\text{G}_{\Q_p}})$. In the literature one has the following results on weight and level lowering.

\begin{theorem}[\cite{BE}, Theorem 4.5]\label{th12} Let $p$ be an odd prime and let $\rho :G_{\Q} \to GL_2(\overline{\F}_p)$ be a continuous, irreducible and odd mod $p$ Galois representation. Suppose that there exists a Katz eigenform $g \in S_k(\Gamma_1(N), \varepsilon, \overline{\mathbb{F}}_p)_{\text{Katz}}$ such that $\rho$ is isomorphic to $\rho_g$. Then there exists a Katz eigenform $f \in S_{k(\rho)}(\Gamma_1(N), \varepsilon, \overline{\mathbb{F}}_p)_{\text{Katz}}$ with the same eigenvalues for $T_l (l \ne p)$ as $g$ has, such that $\rho$ is isomorphic to $\rho_f$. Moreover, there is no eigenform of level prime to $p$ and of weight less than $k(\rho)$ whose associated Galois representation is isomorphic to $\rho$.
\end{theorem}

As a remark, due to the theorem of C. Khare and J.-P. Wintenberger \cite{KW}, \cite{KW1} and of M. Kisin \cite{KM} proving Serre's conjecture there exists a Katz eigenform $F\in S_{k}(\Gamma_1(N), \mathbb{F}_p)_{\text{Katz}}$ for some integers $k$ and $N$ such that $F$ gives rise to the same Galois representation of the above theorem. Thus the existence of $g$ is superfluous. The case $p=2$ is explained in \cite{FC}. For most results we need the following version of Serre's conjecture.

\begin{theorem}[Serre's conjecture, \cite{JS}, \cite{KW}, \cite{KW1} and \cite{KM}]\label{pr132} Let $f\in S_k(\Gamma_1(N), \varepsilon,\\ \overline{\mathbb{F}}_p)_{\text{Katz}}$ be a Katz eigenform such that $\rho_f$ is irreducible. Then there exists a unique Katz newform $g$ in level $N(\rho)$ and weight $k(\rho)$ such that $\rho_f \cong \rho_g$.
\end{theorem}
\begin{proof}
Since $\rho_f$ is irreducible we have by Theorem~\ref{th12} that there exists a cuspidal eigenform $h\in S_{k(\rho)}(\Gamma_1(N), \varepsilon, \overline{\mathbb{F}}_p)_{\text{Katz}}$ such that $\rho_f \cong \rho_h$. Let $g$ be the level lowering of $h$ to level $N(\rho)$ (see \cite{KB} and \cite{KR}). Then by definition of level lowering we have $\rho_h \cong \rho_g$ which gives the result as $f$ is uniquely determined by Theorem~\ref{th02}.
\end{proof}

\begin{prop}\label{pr11} Let $f \in S_k(\Gamma_1(N), \overline{\mathbb{F}}_p)_{\text{Katz}}$ be a Katz modular form such that $f(q) \in \overline{\F}_{p}[[q^l]]$ for some prime number $l\ne p$. Then there exists a unique cusp form $g \in S_k(\Gamma_1(N/l), \overline{\mathbb{F}}_p)_{\text{Katz}}$ such that $f(q)=g(q^l)$. In particular, $g=0$ if $l\nmid N$.
\end{prop}
\begin{proof}
This follows from Lemma~3.6 of \cite{SA} when $l|N$ and the proof of Theorem~1.1 of \cite{OR} when $l\nmid N$. In the latter case, the statement follows from the claim in the 4th paragraph of page 31 of \cite{OR}.
\end{proof}

One has the following weight degeneracy maps on Katz modular forms which are not present in classical modular forms. See [\cite{BG}, \S4, pg. 457] for the details. The first one is multiplying a form by the Hasse invariant which is a Katz modular form $A\in M_{p-1}(\Gamma_1(1), \overline{\F}_p)_{\text{Katz}}$ whose $q$-expansion at the cusp infinity is equal to 1. We denote this map by $\alpha_A$. The other one is taking the Frobenius of a form $f$, $\text{Frob}(f)(q)=f(q^p)$. Multiplying a form by the Hasse invariant does not change the level and the $q$-expansion of the form but adds $p-1$ to the weight. Taking the Frobenius of a form multiplies the weight by $p$ but does not change the level. These two degeneracy maps commute with Hecke operators $T_n$ for all $n$ such that $p\nmid n$.

Let $f\in M_{k}(\Gamma_1(N), \overline{\F}_p)_{\text{Katz}}$ be any Katz modular form. Then to introduce the notion of weight old space corresponding to a form $f$ let us recursively associate a weight to a word formed by the letters $\text{Frob}$ and $A$. Let the empty word have weight $k$. Suppose $m$ is a word of length $n$ and weight $w$. Then set the weight of $A \circ m$ to be $w+p-1$ and the weight of $\text{Frob}\circ m$ to be $pw$. Then the $weight~old~space$ of $f$ in the weight $k'\geq k$ is defined by\\
$\mathcal{W}^{\text{old}}_{k'}(f)=\langle W(f): W\text{ is a word in A and Frob such that } W(f) \text{ has weight }k' \rangle _{\overline{\F}_p}$.

Similarly to the classical modular forms one has the notion of level degeneracy maps on Katz modular forms. Let $f \in M_{k}(\Gamma_1(N), \overline{\F}_p)_{\text{Katz}}$ be any Katz modular form, let $M$ be a positive multiple of $N$ with $\text{gcd}(p,M)=1$ and let $d\geq 1$ be an integer such that $d|M$. Then we have the $d$-th degeneracy map in level $M$,
\begin{equation*}
B_d^{M}: M_{k}(\Gamma_1(M/d), \overline{\F}_p)_{\text{Katz}} \to M_{k}(\Gamma_1(M), \overline{\F}_p)_{\text{Katz}}
\end{equation*}
given by $f(q) \mapsto f(q^d)$. Then the $level~old~space$ of $f$ in the level $M$ is given by
\begin{equation*}
\mathcal{L}^{\text{old}}_{M}(f)=\big\langle \big(B_d^{M} \circ B_1^{M/d}\big)(f): d|M/N \big\rangle _{\overline{\F}_p} \subset M_{k}(\Gamma_1(M), \overline{\F}_p)_{\text{Katz}}.
\end{equation*}

By examining the $q$-expansion, it is clear that we have the following commutativity properties: $B_d^{M}\circ \text{Frob} = \text{Frob} \circ B_d^{M}$ and $B_d^{M}\circ A = A \circ B_d^{M}$. Then the $level~and~weight~old~space$ of $f$ in the level $M$ and weight $k'$ is the $\overline{\F}_p$ vector space generated by $\big(B_d^{M} \circ B_1^{M/d}\circ W\big)(f)$ where $d|M/N$ and $\text{W}$ is a word in $\text{A}$ and $\text{Frob}$ such that $W(f)$ has weight $k'$.

\begin{prop}[\cite{NK}, Corollary 4.4.2]\label{pr13} Let $f\in M_k(\Gamma_1(N), \varepsilon, \overline{\mathbb{F}}_p)_{\text{Katz}}$ and $g \in M_{k'}(\Gamma_1(N), \varepsilon', \overline{\mathbb{F}}_p)_{\text{Katz}}$ be Katz eigenforms with $\rho_f \cong \rho_g$. Then $k\equiv k'(\text{mod}~p-1)$ and $\varepsilon=\varepsilon'$, provided that they are primitive.
\end{prop}

\begin{proof}
Let us set $\tilde{\varepsilon}$ and $\tilde{\varepsilon}'$ to be the corresponding 1-dimensional Galois representations of $\varepsilon$ and $\varepsilon'$. Then since $p\nmid N$, $\tilde{\varepsilon}$ is unramified at $p$ and so $\tilde{\varepsilon}(\text{Frob}_p)=\varepsilon(p)$ is well defined. By restricting $\tilde{\varepsilon}\chi_p^{k-1} = \tilde{\varepsilon}'\chi_p^{k'-1}$ to the inertia group $I_p$ we get $\chi_p^{k-1}=\chi_p^{k'-1}$ from which $k \equiv k'(\text{mod}~p-1)$ follows, so $\tilde{\varepsilon}=\tilde{\varepsilon}'$. We get $\varepsilon=\varepsilon'$, for all primes $l\nmid N$.
\end{proof}

Thus for Katz modular forms $f\in M_k(\Gamma_1(N),\overline{\mathbb{F}}_p)_{\text{Katz}}$ and $f_0 \in M_{k'}(\Gamma_1(N),\overline{\mathbb{F}}_p)_{\text{Katz}}$ with the same $q$-expansions where $k\geq k'$ we have $f=A^tf_0$ where $t=(k-k')/(p-1)$.

There exists a derivation $\theta:M_k(\Gamma_1(N), \overline{\F}_p)_{\text{Katz}}\to M_{k+p+1}(\Gamma_1(N), \overline{\F}_p)_{\text{Katz}}$ which increases weights by $p+1$ and whose effect upon each $q$-expansions is $q\frac{d}{dq}$. For a form $f\in M_k(\Gamma_1(N), \overline{\F}_p)_{\text{Katz}}$, $f$ and $\theta^{p-1}f$ have the same Galois representations. By the principle of $q$-expansions we have that the operator $\theta$ maps modular forms to cusp forms.

\begin{prop}[\cite{NK1}, Corollary~5 and Corollary~6]\label{pr110} Let $f \in M_k(\Gamma_1(N), \overline{\mathbb{F}}_p)_{\text{Katz}}$ be a Katz modular form such that $\theta f=0$. Then we can uniquely write $f(q)=A^r g(q^p)$ with $0\leq r \leq p-1, r+k\equiv 0(\text{mod}~p)$ and $g \in M_l(\Gamma_1(N), \overline{\mathbb{F}}_p)_{\text{Katz}}$ with $pl+r(p-1)=k$. Furthermore, if $f$ is a cusp form, then so is $g$ by its uniqueness.
\end{prop}

\section{Strong multiplicity one}\label{Sec 3}
\setcounter{equation}{0}

Let us start by proving that by moving into higher level we can make some of the inside level coefficients of any Katz eigenform zero. Let $g\in S_k(\Gamma_1(N), \overline{\mathbb{F}}_p)_{\text{Katz}}$ be a cuspidal Katz eigenform. Then $g$ is called an outside $N'$ eigenform if $g$ is eigenform for all $T_n$ where $\text{gcd}(n,N')=1$.

\begin{lemma}\label{le21}
Let $f\in M_k(\Gamma_1(N), \overline{\mathbb{F}}_p)_{\text{Katz}}$ be a normalised Katz eigenform and let $N=\prod_{i=1}^{n}l_i^{\alpha_i}$ be the prime factorization of $N$ with $\alpha_i \geq 1$. Let $I_N=\{l_1, l_2, \dots, l_n \}$ and $S \subset I_N$ be any subset. Then there exists a normalised Katz eigenform $\widetilde{f}\in M_k(\Gamma_1(N\prod_{l_i \in S}l_i), \overline{\mathbb{F}}_p)_{\text{Katz}}$ such that $a_{l}(\widetilde{f})=a_{l}(f)$ for all primes $l\notin S$ and $a_{l^m}(\widetilde{f})=0$ for all $l \in S$ and $m \in \Z_{>0}$.
\end{lemma}
\begin{proof}
Let us define $\widetilde{f}(q):=f(q)-\sum_{l_i \in S} a_{l_i}(f)B_{l_i}^{N\prod_{l_i \in S}l_i}f(q)$. Then by properties of level degeneracy maps $\tilde{f}$ is an outside $\prod_{l_i \in S}l_i$ Katz eigenform. On the other hand from simple $q$-expansion calculations we have $a_{l_i^m}(\widetilde{f})=0$ for all primes $l_i \in S$ and integers $m\geq 1$. Thus $\widetilde{f}$ is an eigenform at primes $l_i \in S$.
\end{proof}

\begin{proof}[Proof of Theorem~\ref{th02}]
Let $\rho_f$ and $\rho_{g}$ be the associated Galois representations. Then by hypothesis, the traces of $\rho_f$ and $\rho_{g}$ agree on the Frobenius elements for all primes in a set of primes of density 1. This implies that $a_l(f)=a_l(g)$ for all primes $l\nmid NN'p$ and $\rho_f \cong \rho_{g}$. Then by definition of cuspidal Katz newforms we have $N=N'$ and $k=k'$. Thus $a_l(f)=a_l(g)$ for all primes $l\nmid pN$. Let $N=\prod_{i=1}^{n}l_i^{\alpha_i}$ be the prime factorization of $N$. Then by taking $S=I_N$ in Lemma~\ref{le21} we have forms $\widetilde{f}$ and $\widetilde{g}$ such that $a_{{l^m_i}}(\widetilde{f})=0=a_{{l^m_i}}(\widetilde{g})$ for all $l_i\in S$ and all $m\geq 1$. If $\widetilde{f}-\widetilde{g} \ne 0$, then by Proposition~\ref{pr110} it must be up to a suitable power multiple of the Hasse invariant in the image of Frobenius of some cusp form of weight smaller than $k$, which is impossible by the minimality of $k=k(\rho)$. Thus, $a_p(f)=a_p(g)$. If $a_{l_i}(f) \ne a_{l_i}(g)$ for some $l_i\in I_N$, then by taking $S=I_N-\{l_i\}$ in above lemma we have cusp forms $\widetilde{f}_i$ and $\widetilde{g}_i$ such that $a_{{l^m_j}}(\widetilde{f}_i) = 0 = a_{{l^m_j}}(\widetilde{g}_i)$ for all $l_j\in S$ and $m\geq 1$. Let $\widetilde{G}:=\widetilde{f}_i-\widetilde{g}_i$. Then by Proposition~\ref{pr11}, $\widetilde{G}(q)=\widetilde{G}_1(q^{l_i})$ for some cusp form $\widetilde{G}_1$ of level $\frac{N}{l_i}\prod_{l_j \in S}l_j$ which is impossible by the $l_i$-minimality. Thus $f=g$. Furthermore, for all $d\in \N$ such that $\text{gcd}(d,N)=1$ we have $\langle d \rangle f= \varepsilon(d)\cdot f$ and $\langle d \rangle g= \varepsilon'(d)\cdot g$. Then $f=g$ gives $\varepsilon=\varepsilon'$.
\end{proof}

\begin{lemma}\label{le2}(i). Let $f \in S_k(\Gamma_1(N), \overline{\mathbb{F}}_p)_{\text{Katz}}$ be a normalised Katz eigenform with $p_1$-minimal level and $g \in S_{k}(\Gamma_1(Np_1^{m_1}), \overline{\mathbb{F}}_p)_{\text{Katz}}$ where $m_1\geq 1$ be an outside $p_1$ eigenform such that for all positive integers $n$ such that $p_1\nmid n$, $T_ng=a_n(f)g$. Then $g\in \mathcal{L}^{\text{old}}_{Np_1^{m_1}}(f)$.\\
(ii). Let $M/N=\prod_{i=1}^{t}p_i^{m_i}$ be the prime factorization of $M/N$ where $m_i\geq 1$. Let $f\in S_{k}(\Gamma_1(N), \overline{\mathbb{F}}_p)_{\text{Katz}}$ be a normalised Katz eigenform with $p_1, p_2, p_3, \dots, p_t$-minimal level and $g\in S_{k}(\Gamma_1(M), \overline{\mathbb{F}}_p)_{\text{Katz}}$ be a normalised Katz eigenform satisfying $T_l(g)=a_l(f)g$ for all primes $l\ne p_1, p_2, p_3, \dots, p_t$. Then $g\in \mathcal{L}^{\text{old}}_{M}(f)$.
\end{lemma}

\begin{proof}
(i). We have $a_n(g-a_1(g)f)=0$ for all integers $n\geq 1$ such that $p_1\nmid n$ because $a_n(g-a_1(g)f)=a_1(T_n(g-a_1(g)f))=0$. Then by Proposition~\ref{pr11}, $(g-a_1(g)f)(q)=F(q^{p_1})$ for some outside $p_1$ eigenform $F$ of level $Np_1^{m_1-1}$ such that $\rho_g\cong \rho_F$. Then we proceed by induction on $m_1$. When $m_1=0$ by $p_1$-minimality we have $F=0$ and $g=a_1(g)f$. Assume the result holds for $m_1$ less than some positive integer $m$. Then when $m_1=m$ we have $(g-a_1(g)f)(q)=G(q^{p_1})$ for some outside $p_1$ eigenform $G$ of level $Np_1^{m-1}$ such that $\rho_g\cong \rho_G$, which by induction assumption is in the level old space of $f$ in level $Np_1^{m-1}$. Thus $g\in \mathcal{L}^{\text{old}}_{Np_1^{m_1}}(f)$.\\
(ii). The case $t=1$ follows from (i) above. Assume the result holds for $t<r$. Let $f \in S_{k}(\Gamma_1(N), \overline{\mathbb{F}}_p)_{\text{Katz}}$ be a normalised eigenform with $p_1, p_2, p_3, \dots, p_r$-minimal level and $g \in S_{k}(\Gamma_1(M), \overline{\mathbb{F}}_p)_{\text{Katz}}$ be a normalised eigenform with $T_lg=a_l(f)g$ for all primes $l\ne p_1, p_2, p_3, \dots, p_r$. Then by using Lemma~\ref{le21} with $S=\{p_1\}$ we have normalised eigenforms $\widetilde{f} \in S_{k}(\Gamma_1(Np_1^2), \overline{\mathbb{F}}_p)_{\text{Katz}} \subset S_{k}(\Gamma_1(Np_1^{m_1+2}), \overline{\mathbb{F}}_p)_{\text{Katz}}$ and $\widetilde{g}\in S_{k}(\Gamma_1(Mp_1^2), \overline{\mathbb{F}}_p)_{\text{Katz}}$ such that $\widetilde{f}$ has $p_2, p_3, \dots, p_r$-minimal level and $\widetilde{g}$ satisfies $T_l\widetilde{g}=a_l(\widetilde{f})\widetilde{g}$ for all primes $l\ne p_2, p_3, \dots, p_r$. Then by the induction assumption $\widetilde{g}\in \mathcal{L}^{\text{old}}_{Mp_1^2}(\widetilde{f})$, say $\widetilde{g}(q)=\sum_{d|M/(Np_1^{m_1})} \beta_d \widetilde{f}(q^d)$. Here we have assumed that $\widetilde{f}\in S_{k}(\Gamma_1(Np_1^{m_1+2}), \overline{\mathbb{F}}_p)_{\text{Katz}}$. Let $h(q):=\sum_{d|M/(Np_1^{m_1})}\beta_d f(q^d) \in S_{k}(\Gamma_1(Np_2^{m_2}\cdots p_r^{m_r}), \overline{\mathbb{F}}_p)_{\text{Katz}}$. Then $h$ is a normalised outside $M/Np_1^{m_1}$ eigenform with $p_1$-minimal level. On the other hand, since $a_{l^n}(h)=a_{l^n}(\widetilde{g})=a_{l^n}(g)$ for any positive integer $n\geq 1$ and primes $l= p_2, p_3, \dots, p_r$ and $\widetilde{g}$ is an eigenform we have by simple $q$-expansion calculations that $h$ is an eigenform at primes $l=p_2, p_3, \dots, p_r$. Thus $h$ is a normalised Katz eigenform with $p_1$-minimal level. Then $T_lg=a_l(h)g$ for all primes $l\ne p_1$. Then by part (i) above we have $g\in \mathcal{L}^{\text{old}}_{M}(h)$, so $g \in \mathcal{L}^{\text{old}}_{M}(f)$.
\end{proof}

\begin{lemma}\label{co24} Let $f \in S_{k}(\Gamma_1(N, \overline{\mathbb{F}}_p)_{\text{Katz}}$ be a Katz newform. Then any normalised Katz eigenform $g \in S_{k}(\Gamma_1(M), \overline{\mathbb{F}}_p)_{\text{Katz}}$ such that $\rho_f \cong \rho_{g}$ is in the level old space of $f$.
\end{lemma}

\begin{proof}
By the hypothesis we have $f=B^{M}_{1}f \in S_k(\Gamma_1(M), \overline{\mathbb{F}}_p)_{\text{Katz}}$ as $N|M$. Then setting $S=I_M$ in Lemma~\ref{le21} and applying Proposition~\ref{pr110} we have $(\widetilde{f}-\widetilde{g})(q)=A^rG(q^p)$ for some integer $r$ and an outside $p$ Katz eigenform $G$ of weight smaller than $k$, which is impossible by the minimality of weight unless $G=0$, so we have $a_p(f)=a_p(g)$. Suppose $a_l(f)\ne a_l(g)$ for some prime $l\nmid M/N$ and $l|M$. Then taking $S=I_M-\{l\}$ in Lemma~\ref{le21} gives forms $\widetilde{f}$ and $\widetilde{g}$ such that $a_{l'}(\widetilde{f}) = 0 = a_{l'}(\widetilde{g})$ for all primes $l'|M$ and $l'\ne l$. Then $\widetilde{f}(q)- \widetilde{g}(q)=F(q^l)$ for some modular form $F\ne 0$ of level $\frac{M}{l}\prod_{l' \in S}l'$ which is impossible by $l$-minimality. Thus $T_lg=a_l(f)g$ for all primes $l\nmid M/N$. Then since the level $N$ of $f$ is $l$-minimal for any prime $l$, in particular it is $l$-minimal for $l| M/N$. Then by applying Lemma~\ref{le2} we have $g \in \mathcal{L}^{\text{old}}_{M}(f)$.
\end{proof}

\begin{cor}\label{pr25}
Let $f\in S_k(\Gamma_1(N), \overline{\mathbb{F}}_p)_{\text{Katz}}$ be a normalised Katz eigenform with $p_1$-minimal level for all $p_1|N$. Then any normalised Katz eigenform $g\in S_k(\Gamma_1(M),\\ \overline{\mathbb{F}}_p)_{\text{Katz}}$ such that $\rho_f \cong \rho_g$ and $a_p(f)=a_p(g)$ is in the level old space of $f$.
\end{cor}

Let $f \in M_{k}(\Gamma_1(N), \overline{\mathbb{F}}_p)_{\text{Katz}}$ be a normalised Katz eigenform with minimal weight $k$. Then for $k'\geq k$, define $V_{f,k'}$ as the $\overline{\mathbb{F}}_p$ vector space generated by $F\in M_{k'}(\Gamma_1(N),\\ \overline{\mathbb{F}}_p)_{\text{Katz}}$ such that $F$ is an outside $p$ Katz eigenform with eigenvalues $\lambda_l(F)=a_l(f)$ for all primes $l\ne p$. Then we have the following
\begin{lemma}\label{le1}
The space $V_{f,k'}$ is a subspace of the weight old space of $f$ in the weight $k'$.
\end{lemma}
\begin{proof}
We proceed by induction on $k'$. Let $k'=k$. Then for every $F \in V_{f,k'}$, by Proposition~\ref{pr110}, we can write $(F-a_1(F)f)(q)=A^rG(q^p)$ for some integer $r$ and an outside $p$ Katz eigenform $G$ of weight smaller than $k'$, which is impossible by the minimality of weight unless $G=0$, so we have $V_{f,k'} = \langle f \rangle$. Then suppose the induction hypothesis is correct for all weights less than $k'$. Then by Proposition~\ref{pr13}, $k'=k+m(p-1)$ for some non-negative integer $m$. Set $f_0=A^mf\in V_{f,k'}$. Then since $V_{f,k'}$ is a finite dimensional $\overline{\mathbb{F}}_p$-vector space, say of dimension $d$, we can pick modular forms $f_1, f_2, f_3, \dots, f_{d-1}\in V_{f,k'}$ such that $f_0, f_1, f_2, \dots, f_{d-1}$ constitutes a basis for $V_{f,k'}$. Then for all $1\leq i\leq d-1$, define $g_i:=f_i-a_1(f_i)f_0$. Then $a_1(g_i)=0$ which gives $a_n(g_i)=0$ for all integers $n\geq 1$ such that $p\nmid n$ as $a_n(g_i)=a_1(T_n g_i)=a_n(f)a_1(g_i)=0$. Then by Proposition~\ref{pr110} there exist modular forms $\widetilde{g}_i\in M_{k_i}(\Gamma_1(N), \overline{\mathbb{F}}_p)_{\text{Katz}}$ for $i=1, 2, 3, \dots, d-1$ such that $g_i(q)=A^{r_i}\widetilde{g}_i(q^p)$ and $\widetilde{g}_i \in V_{f,k_i}$ for some integers $r_i$ and $k_i< k'$. Then by the induction assumption $\widetilde{g}_1, \widetilde{g}_2, \widetilde{g}_3, \dots, \widetilde{g}_{d-1}$ are in the weight old space of $f$ in weights $k_1, k_2, k_3, \dots, k_{d-1}$ respectively. This implies that the basis elements $f_1, f_2, f_3, \dots, f_{d-1}$ are in the weight old space of $f$ in weight $k'$. This gives the result.
\end{proof}

\begin{cor}\label{le25} Let $f \in M_{k}(\Gamma_1(N), \overline{\mathbb{F}}_p)_{\text{Katz}}$ be a normalised Katz eigenform with minimal weight. Then any normalised Katz eigenform $g\in M_{k'}(\Gamma_1(N), \overline{\mathbb{F}}_p)_{\text{Katz}}$ such that $\rho_f\cong \rho_g$ and $a_l(f)=a_l(g)$ for all primes $l|N$ is in the weight old space of $f$.
\end{cor}

In the above corollary one cannot relax the condition that the eigenvalues $a_l(f)$ and $a_l(g)$ for $T_l$ for all primes $l$ dividing the level are the same. To construct a counterexample let $F \in S_{k}(\Gamma_1(M), \overline{\mathbb{F}}_p)_{\text{Katz}}$ be a Katz newform. Then choose a prime $l\nmid Mp$ such that $T_l$ has two distinct eigenvalues on $\langle F(q), F(q^l) \rangle \subseteq S_{k}(\Gamma_1(Ml), \overline{\mathbb{F}}_p)_{\text{Katz}}$. Then we can produce normalised Katz eigenforms $f$ and $g$ in this subspace such that $\rho_f \cong \rho_g$ and $a_l(f)\ne a_l(g)$ but one is not in the weight old space of the other.

In Proposition~\ref{pr132} we have associated a Katz newform to any Katz eigenform which has an irreducible Galois representation. More generally if we assume the existence of Katz newforms for Katz eigenforms which has reducible Galois representation we have

\begin{proof}[Proof of Theorem~\ref{th01}] By applying Lemma~\ref{le21} with $S=I_M$ we have eigenforms $\widetilde{F}$ and $\widetilde{f}$ such that $a_l(\widetilde{F})=0=a_l(\widetilde{f})$ for all primes $l|M$. Then by using Corollary~\ref{le25} one can write $\widetilde{F}(q)=\sum_{\delta \in D_{k}^{k'}}\alpha_{\delta} \delta(\widetilde{f}(q))$ for some $\alpha_{\delta} \in \overline{\mathbb{F}}_p$ where $D_{k}^{k'}$ is the set of words $W$ in $A$ and $\text{Frob}$ such that $W$ takes weight $k$ forms into weight $k'$ forms. Let us define $F_1(q):=\sum_{\delta \in D_{k}^{k'}}\alpha_{\delta} \delta(f(q))$ by replacing the form $\widetilde{f}$ by $f$. Suppose that $f$ has Dirichlet character $\varepsilon$ and suppose that $F_1(q):=\sum_{\substack{t=0 \\ i}}\beta_t A^i f(q^{p^t})$. Then $F_1$ is a $T_p$ Katz eigenform since we have $a_{pm}(F_1)=a_p(F_1)a_m(F_1)$ for any positive integer $m$ such that $\text{gcd}(m, p)=1$ and $a_{p^n}(F_1)=a_p(F_1)a_{p^{n-1}}(F_1)- p^{k'-1}\varepsilon(p)a_{p^{n-2}}(F_1)$ for any positive integer $n\geq 2$. The last relation follows by a case by case calculation using the later definition of $F_1$. Then since $\rho_{F}\cong \rho_{F_1}, a_{p}(F)=a_{p}(F_1)$ and level$(F_1)=N$, by applying Corollary~\ref{pr25} we have $F \in \mathcal{L}^{\text{old}}_{M}(F_1)$, which gives the desired result.
\end{proof}

In particular, we have determined all possible coefficients of any Katz eigenform which has irreducible mod $p$ Galois representation from the coefficients of the corresponding Katz newform.

\begin{cor}\label{le26}
Let $F\in S_{k}(\Gamma_1(M), \overline{\mathbb{F}}_p)_{\text{Katz}}$ be a normalised Katz eigenform and assume that there exists a Katz newform $f\in S_{k'}(\Gamma_1(N), \varepsilon, \overline{\mathbb{F}}_p)_{\text{Katz}}$ such that $\rho_F \cong \rho_f$. Suppose that $F(q)=\sum_{n\geq 1}a_nq^n$ and $f(q)=\sum_{n\geq 1}b_nq^n$ are their $q$-expansions. Then we have the following identities:\\
(i). When prime $l\nmid Mp/N$, $a_l=b_l$.\\
(ii). When prime $l|M/N$ and $l|N$, $a_l=0$ or $a_l=b_l$.\\
(iii). When prime $l|Mp/N$ but $l\nmid N$, $a_l=0$ or $a_l^2-a_lb_l+\varepsilon(l)l^{k-1}=0$.
\end{cor}
\begin{proof}
The result follows from the explicit $q$-expansion calculations of the relations of Theorem~\ref{th01}.
\end{proof}
\begin{rem}
Let us recall that in classical newform theory, the newspace has a basis consisting of newforms. However this cannot be generalised to Katz modular forms. A counterexample occurs in $S_{1}(\Gamma_0(229), \overline{\mathbb{F}}_2)_{\text{Katz}}$. The associated Hecke algebra $\mathbb{T}$ is a local 2-dimensional $\overline{\mathbb{F}}_2$-algebra, hence it has a unique attached Katz eigenform, whereas $S_{1}(\Gamma_0(229), \overline{\mathbb{F}}_2)_{\text{Katz}}$ is a 2-dimensional $\overline{\mathbb{F}}_2$-algebra, which is a contradiction.
\end{rem}
\section{Reducible case}\label{Sec 4}
\setcounter{equation}{0}

In this section, we will prove a strong multiplicity one result for Katz Eisenstein series. Then later we will show that under some condition a reducible mod $p$ Galois representation arises from a normalised Katz eigenform with optimal level.

Let us start by defining and studying the properties of Eisenstein series. We will follow the notations and definitions of [\cite{WS}, Chapter 5].

We define the generalized Bernoulli number $B_k^{\varepsilon}$ attached to a complex modulo $n$ Dirichlet character $\varepsilon$ by the following infinite series
\begin{equation*}
\sum_{j=1}^n\frac{\varepsilon(j)xe^{jx}}{e^{nx}-1}=\sum_{k=0}^{\infty}\frac{B_k^{\varepsilon}x^k}{k!}.
\end{equation*}

If $\varepsilon$ is the trivial character, then $B_k^{\varepsilon}$ is equal to the classical Bernoulli number $B_k$ for $k>1$ and $B_1^{\varepsilon}=-B_1=1/2$.

Let $\varepsilon_1$ and $\varepsilon_2$ be two Dirichlet characters modulo $u$ and $v$ such that $uv|N$ and let $k$ be a positive integer such that $(\varepsilon_1 \varepsilon_2)(-1) = (-1)^k$. Let $t$ be a positive integer. Then we have the Eisenstein series $E_k^{\varepsilon_1,\varepsilon_2}(q)$ defined by a power series
\begin{equation*}
E_k^{\varepsilon_1,\varepsilon_2}(q) := c_0 + \sum_{m\geq 1} \bigg( \sum_{0<d|m}\varepsilon_1(d)\varepsilon_2\bigg(\frac{m}{d}\bigg)d^{k-1}\bigg)q^m \in \mathbb{Q} (\varepsilon_1,\varepsilon_2)[[q]]
\end{equation*}
where $c_0 =-\frac{B_k^{\varepsilon_1}}{2k}$ when $\text{cond}(\varepsilon_2) = 1$ and $c_0 =0$ otherwise. Then, except when $k = 2$ and $\varepsilon_1 = \varepsilon_2 = 1$, the power series $E_k^{\varepsilon_1,\varepsilon_2}(q^t)$ belongs to $M_k(\Gamma_1(tuv), \mathbb{C})$ for all $t\geq 1$. If $k = 2$ and $\varepsilon_1 = \varepsilon_2 = 1$, let $t >1$, then $E_2^{1,1}(q) - tE_2^{1,1}(q^t)$ is a modular form in $M_2(\Gamma_1(t), \mathbb{C})$. Moreover the modular form $E_k^{\varepsilon_1,\varepsilon_2}(q)$ is a normalized eigenform for all Hecke operators. Analogously, for all positive integers $t > 1$ the series $E_2^{1,1}(q) - tE_2^{1,1}(q^t)$ is a normalised eigenform for all Hecke operators. Let us set $E_k^{\varepsilon_1,\varepsilon_2, t}(q)$ to $E_2^{\varepsilon_1,\varepsilon_2}(q) - tE_2^{\varepsilon_1,\varepsilon_2}(q^t)$ when $k=2$ and $\varepsilon_1=\varepsilon_2=1$, and to $E_k^{\varepsilon_1,\varepsilon_2}(q^t)$ otherwise.

Hereafter we will assume that all Dirichlet characters which we consider are primitive and we are not in the situation where our Eisenstein series have weight $k=2$ and level $N=1$ since there is no holomorphic modular form of such parameters. Let $\text{Den}\Big(\frac{B_k^{\varepsilon_1}}{2k}\Big)$ denotes the denominator of $\frac{B_k^{\varepsilon_1}}{2k}$ when it is written in a reduced fractional form.

We can define Katz Eisenstein series by considering the mod $p$ reduction of Eisenstein series. For prime $p\nmid \text{Den}\Big(\frac{B_k^{\varepsilon_1}}{2k}\Big)$ when $\text{cond}(\varepsilon_2)=1$, we define the Katz eigenform $E_k^{\varepsilon, \varepsilon'}$ as the mod $p$ reduction of the associated Eisenstein series $E_k^{\varepsilon_1, \varepsilon_2}$. More precisely, we may first assume that the level of the classical Eisenstein series $N$ is at least 5, since for some fixed prime $l\nmid 6Np$ we have the level degeneracy map $B_1^{lN}$ which increases the level by a multiple of $l$. Then we can observe that all coefficients of $B_1^{lN}E_k^{\varepsilon_1, \varepsilon_2}$ belongs to $\mathcal{O}$ where $\mathcal{O}$ is the ring of integer of some finite extension of $\Q_p$. Then by [{\cite{DI}, Theorem 12.3.4.2}] we can conclude that $E_k^{\varepsilon_1, \varepsilon_2}$ is a modular form with coefficients in $\mathcal{O}$, i.e., $E_k^{\varepsilon_1, \varepsilon_2}$ has $p$-integral coefficients at all cusps. Then by the condition ensuring $p\nmid \text{Den}\Big(\frac{B_k^{\varepsilon_1}}{2k}\Big)$ when $\text{cond}(\varepsilon_2)=1$, we have a well defined reduction of $E_k^{\varepsilon_1, \varepsilon_2}$ at some prime above $p$, $E_k^{\varepsilon, \varepsilon'}$ which is a Katz Eisenstein series. Here $\varepsilon$ and $\varepsilon'$ are the mod $p$ reductions of $\varepsilon_1$ and $\varepsilon_2$ when considered as Galois representations. Similarly by taking different positive integers $t \geq 1$ we can define a Katz Eisenstein series $E_k^{\varepsilon, \varepsilon', t}$ as a mod $p$ reduction of $E_k^{\varepsilon_1, \varepsilon_2, t}$. A normalized Katz eigenform $f(q)=E_k^{\varepsilon, \varepsilon'}(q) \in M_k(\Gamma_1(N), \overline{\mathbb{F}}_p)_{\text{Katz}}$ is called a \emph{Katz~new~Eisenstein~series} if it satisfies the condition of Definition~\ref{df01}.

We know that a reducible mod $p$ Galois representation arises from Eisenstein series. Let $f=E_k^{\varepsilon, \varepsilon', t}\in M_k(\Gamma_1(N), \chi, \overline{\mathbb{F}}_p)_{\text{Katz}}$ be any Katz eigenform. Then we have $\rho_f\cong \varepsilon' \oplus \varepsilon \chi_p^{k-1}$ which is unramified outside $pN$ with the property that $\text{tr}(\rho_f(\text{Frob}_l))=f(T_l)$ and $\det (\rho_f(\text{Frob}_l))=\chi(l)l^{k-1}$ for all primes $l\nmid pN$. In fact the converse also holds. Any reducible mod $p$ Galois representation comes from some twist of an Eisenstein series.

Let $\varepsilon$ and $\varepsilon'$ be primitive Dirichlet characters with values in $\overline{\mathbb{F}}_p$ and let $\varepsilon_1$ and $\varepsilon_2$ be their respective complex lifings with the same conductors and the same orders. Then we start by proving the existence of Katz new Eisenstein series.

\begin{prop}\label{pr31}
Let $1\leq k\leq p-1$ and $N=\text{cond}(\varepsilon)\cdot\text{cond}(\varepsilon')$ be positive integers. Assume $k\ne 2$ if $\varepsilon=\varepsilon'=1$ and $p\nmid \text{Den}\Big(\frac{B_k^{\varepsilon_1}}{2k}\Big)$ if $\varepsilon_2=1$. Then $E_k^{\varepsilon, \varepsilon'}\in M_k(\Gamma_1(N), \overline{\mathbb{F}}_p)_{\text{Katz}}$ is a Katz new Eisenstein series such that $\rho_f\cong \varepsilon' \oplus \varepsilon \chi_p^{k-1}$.
\end{prop}
\begin{proof}
This immediately follows from the discussion above. Here $E_k^{\varepsilon, \varepsilon'}$ is Katz new Eisenstein series as it a normalised eigenform with optimal level and weight.
\end{proof}

\begin{cor}
Let $F\in M_k(\Gamma_1(N), \overline{\mathbb{F}}_p)_{\text{Katz}}$ be a normalised Katz eigenform such that $\rho_F\cong \varepsilon' \oplus \varepsilon \chi_p^b$ where $0\leq b\leq p-2$. Suppose $N=\text{cond}(\rho_F)$. Assume $b\not\equiv 1(\text{mod}~p-1)$ if $\text{cond}(\rho_F)=1$ and $p\nmid \text{Den}\Big(\frac{B_{b+1}^{\varepsilon_1}}{2(b+1)}\Big)$ if $\varepsilon_2=1$. Then $F$ is in the weight old space of $E_{b+1}^{\varepsilon, \varepsilon'}$ in the weight $k$.
\end{cor}
\begin{proof}
By Proposition~\ref{pr31} above, $f=E_{b+1}^{\varepsilon, \varepsilon'}$ is a Katz new Eisenstein series such that $\rho_f \cong \rho_F$. By Proposition~\ref{pr13} we can write $k=(b+1)+m(p-1)$ for some non-negative integer $m$. Then by comparing coefficients of $\theta^{p-1}F$ and $\theta^{p-1}A^mf$ using Corollary~\ref{pr25} we have $a_l(F)=a_l(f)$ for all primes $l|N$. Then Corollary~\ref{le25} completes the proof.
\end{proof}

\begin{prop}\label{th32}
Let $f(q)=E_k^{\varepsilon, \varepsilon'}(q) \in M_k(\Gamma_1(N), \overline{\mathbb{F}}_p)_{\text{Katz}}$ and $g(q)=E_{k'}^{\chi, \chi'}(q) \in M_{k'}(\Gamma_1(N'), \overline{\mathbb{F}}_p)_{\text{Katz}}$ be Katz new Eisenstein series with $a_l(f)=a_l(g)$ for each $l$ in a set of primes of density 1. Then $f=g, k=k', N=N', \varepsilon = \chi$ and $\varepsilon' = \chi'$ or $\varepsilon' = \chi$ and $\varepsilon = \chi'$ when $k\equiv 1(\text{mod}~p-1)$.
\end{prop}
\begin{proof}
From the hypothesis we have $\rho_{f} \cong \rho_{g}$ or $\varepsilon' \oplus \varepsilon \chi_p^{k-1} \cong \chi' \oplus \chi \chi_p^{k'-1}$. Then by the definition of Katz new Eisenstein series the levels and weights of the forms are optimal so they are equal. Then from $\varepsilon' \oplus \varepsilon \chi_p^{k-1} \cong \chi' \oplus \chi \chi_p^{k-1}$ we have the following two cases. (i). $\varepsilon'=\chi'$ and $\varepsilon \chi_p^{k-1} = \chi \chi_p^{k-1}$ or (ii). $\varepsilon' = \chi \chi_p^{k-1}$ and $\varepsilon \chi_p^{k-1} = \chi'$. The first case gives $\varepsilon'=\chi'$ and $\varepsilon=\chi$ as a Galois representations while the second case gives $\varepsilon' = \chi$ and $\varepsilon = \chi'$ provided that $k\equiv 1,p(\text{mod}~p-1)$. This completes the proof.
\end{proof}

\begin{rem}\label{rk1}
It is not always the case to obtain a cuspidal eigenform with both optimal weight and optimal level which gives rise to a given reducible mod $p$ Galois representation. For example, there exists a modular form $f\in S_{28}(\Gamma_1(1), \overline{\mathbb{F}}_7)_{\text{Katz}}$ such that $\rho_f \cong \rho_{\overline{E}_4}$ but there is no cuspidal eigenform $g\in S_{4}(\Gamma_1(1), \overline{\mathbb{F}}_7)_{\text{Katz}}$ such that $\rho_g \cong \rho_f$. Thus for cuspidal eigenforms which have reducible mod $p$ Galois representation we consider Katz Eisenstein series as a candidate of an associated Katz newforms.
\end{rem}

Now let us consider the case when the mod $p$ representation of the modular form in Theorem~\ref{th01} is reducible. To be precise let $F\in S_{k'}(\Gamma_1(M), \overline{\mathbb{F}}_p)_{\text{Katz}}$ be a normalised Katz eigenform with reducible mod $p$ Galois representation $\rho_F \cong \varepsilon' \chi_p^a \oplus \varepsilon \chi_p^b$ where $\det \rho_F=\varepsilon \varepsilon' \chi_p^{k-1}$ for some $a,b\in \mathbb{Z}/(p-1)\mathbb{Z}$ such that $1 \leq k \leq p+1$ and $k-1\equiv b-a(\text{mod}~p-1)$. Assume that $\varepsilon$ and $\varepsilon'$ are primitive when considered as Dirichlet characters.

Then we show that there exists a normalised Katz eigenform $g$ of optimal level such that $\rho_F \cong \rho_g$.

\begin{lemma}\label{le17}
Let $F\in S_{k'}(\Gamma_1(M), \overline{\mathbb{F}}_p)_{\text{Katz}}$ be a normalised Katz eigenform such that $\rho_F \cong \varepsilon'\chi_p^a \oplus \varepsilon \chi_p^b$ where $1 \leq a \leq p-1, 0 \leq b \leq p-2$ and $k$ is a positive integer such that $1 \leq k \leq p+1$. Assume that $(N(\rho_F), k)\ne (1,2)$. Then we have the following four cases:\\
(i). If $\varepsilon_2\ne 1$ and $k-1\equiv b-a(\text{mod}~p-1)$, then $\rho_F$ comes from $g=\theta^a(\overline{E_k^{\varepsilon_1,\varepsilon_2}})$.\\
(ii). If $\varepsilon = \varepsilon' =1, k=2, b-a\equiv 1(\text{mod}~p-1)$ and $p\ne 2, 3$, then $\rho_F$ comes from $g=\theta^{a}(\overline{E_{p^2+1}^{1, 1}})$.\\
(iii). If $\varepsilon_1=\varepsilon_2=1, k-1\equiv b-a(\text{mod}~p-1)$, $k\ne 2$ and $p\nmid \text{Den}\Big(\frac{B_k}{2k}\Big)$, then $\rho_F$ comes from $g=\theta^{a}(\overline{E_k^{1, 1}})$.\\
(iv). If $\varepsilon_1\ne 1, \varepsilon_2 = 1, k-1\equiv b-a(\text{mod}~p-1)$ and $p\nmid \text{Den}\Big(\frac{B_k^{\varepsilon_1}}{2k}\Big)$, then $\rho_F$ comes from $g= \theta^a(\overline{E_k^{\varepsilon_1, 1}})$.
\end{lemma}
\begin{proof}
In the case when $a=p-1$ we have the following cases. (i). We have that $\rho_F$ comes from $E_k^{\varepsilon,\varepsilon', t}$ for some positive integer $t$. Then since $(N(\rho),k)\ne (1,2)$ taking $t=1$ gives a normalised eigenform $E_k^{\varepsilon,\varepsilon'}$ with an optimal level. Here $a_0(E_k^{\varepsilon_1,\varepsilon_2})=0$, so we can take modulo $p$ reduction and apply the theta operator to get a normalised cuspidal Katz eigenform $g=\theta^{p-1}(\overline{E_k^{\varepsilon_1,\varepsilon_2}})$ such that $\rho_F\cong\rho_g$.\\
(ii). Let $\rho_F\cong \overline{\rho}_{E_2^{1,1,t}}$ for some positive integer $t$. Then for prime $p\ne 2,3$ we have $\rho_{\theta^{p-1}\big(\overline{E_2^{1,1,t}}\big)}\cong \rho_{\overline{E_{p^2+1}^{1,1}}}$. Then set $g=\theta^{p-1}\big(\overline{E_{p^2+1}^{1,1}}\big)$.\\
(iii). Here the assumption $p\nmid \text{Den}\Big(\frac{B_k^{\varepsilon_1}}{2k}\Big)$ implies that the modulo $p$ reduction is well defined. Similarly (iv) holds.\\
On the other hand, when $a\ne p-1$, by applying the above method to the twist $\theta^{-a}F$ one can get the remaining results.
\end{proof}

Let us assume that we are in the same notation and under the same assumptions as in Lemma~\ref{le17}. Then as an immediate consequence of Corollary~\ref{pr25} we have

\begin{theorem}\label{th47}
Let $F\in S_{k'}(\Gamma_1(M), \overline{\mathbb{F}}_p)_{\text{Katz}}$ be the above normalised Katz eigenform which we consider. Suppose that $a_p(F)=0$ and $g \in S_k(\Gamma_1(N(\rho_F)), \overline{\mathbb{F}}_p)_{\text{Katz}}$ is the modular form associated to $F$ as in Lemma~\ref{le17}. Then up to a suitable power multiple of the Hasse invariant, $F$ is in the level old space of $g$.
\end{theorem}

Here also as application of Corollary~\ref{le26} we can compute all possible coefficients of any cuspidal Katz eigenform which has reducible mod $p$ Galois representation in terms of the associated normalised Katz eigenform of optimal level.

\begin{cor}
Let $F(q)=\sum_{n\geq 1}a_nq^n \in S_{k}(\Gamma_1(M), \overline{\mathbb{F}}_p)_{\text{Katz}}$ be an ordinary normalised Katz eigenform such that $2\leq k\leq p$ and $\rho_F$ is tamely ramified at $p$ and let $G(q)=\sum_{n\geq 1}b_nq^n \in S_{p+1-k}(\Gamma_1(M), \overline{\mathbb{F}}_p)_{\text{Katz}}$ be the corresponding companion form. Suppose that $f(q)=\sum_{n\ge 0}c_nq^n \in M_{k'}(\Gamma_1(N), \chi, \overline{\mathbb{F}}_p)_{\text{Katz}}$ is given by $\overline{E_{k'}^{\varepsilon_1,\varepsilon_2}}$ when $\rho_F \cong \varepsilon' \oplus \varepsilon \chi_p^{k'-1}$ and by the corresponding Katz newform when $\rho_F$ is irreducible. Then we have the following identities:\\
(i). When prime $l\nmid Mp/N$, $l^kb_l=lc_l$.\\
(ii). When prime $l|M/N$ and $l|N$, $l^kb_l=0$ or $l^kb_l=lc_l$.\\
(iii). When prime $l|M/N$ but $l\nmid N$, $l^kb_l=0$ or $lb_l(l^{k-1}b_l-c_l)+\chi(l)l=0$.
\end{cor}
\begin{proof}
The existence of $G$ is a standard result, see \cite{BG} and \cite{CV}. The identities follows by applying Corollary~\ref{le26} to $n^kb_n=na_n$ for all $n\geq 1$. When $\rho_F$ is reducible we instead compare $G$ with $\theta^{p-1}f$.
\end{proof}

\subsection*{Acknowledgments}
\setcounter{equation}{0}
The author would like to thank Gabor Wiese for his valuable remarks which improves the paper, and Ian Kiming and Samuele Anni for their comments on the manuscript.

\bibliographystyle{abbrv}
\bibliography{Paper}

\begin{thebibliography}{10}

\bibitem{SA}
S.~Anni.
\newblock A note on the minimal level of realization for a mod {$\ell$}
  eigenvalue system.
\newblock In {\em Automorphic forms and related topics}, volume 732 of {\em
  Contemp. Math.}, pages 1--13. Amer. Math. Soc., Providence, RI, 2019.

\bibitem{AL}
A.~O.~L. Atkin and J.~Lehner.
\newblock Hecke operators on {$\Gamma _{0}(m)$}.
\newblock {\em Math. Ann.}, 185:134--160, 1970.

\bibitem{KB}
K.~Buzzard.
\newblock On level-lowering for mod 2 representations.
\newblock {\em Math. Res. Lett.}, 7(1):95--110, 2000.

\bibitem{FC}
F.~Calegari.
\newblock Is serre’s conjecture still open?
\newblock In {\em Blogpost}.
  https://galoisrepresentations.wordpress.com/2014/08/10/is-serres-conjecture-still-open.

\bibitem{CV}
R.~F. Coleman and J.~F. Voloch.
\newblock Companion forms and {K}odaira-{S}pencer theory.
\newblock {\em Invent. Math.}, 110(2):263--281, 1992.

\bibitem{DS}
P.~Deligne and J.-P. Serre.
\newblock Formes modulaires de poids {$1$}.
\newblock {\em Ann. Sci. \'{E}cole Norm. Sup. (4)}, 7:507--530 (1975), 1974.

\bibitem{DI}
F.~Diamond and J.~Im.
\newblock Modular forms and modular curves.
\newblock In {\em Seminar on {F}ermat's {L}ast {T}heorem ({T}oronto, {ON},
  1993--1994)}, volume~17 of {\em CMS Conf. Proc.}, pages 39--133. Amer. Math.
  Soc., Providence, RI, 1995.

\bibitem{BE}
B.~Edixhoven.
\newblock The weight in {S}erre's conjectures on modular forms.
\newblock {\em Invent. Math.}, 109(3):563--594, 1992.

\bibitem{Bj}
B.~Edixhoven and J.-M. Couveignes, editors.
\newblock {\em Computational aspects of modular forms and {G}alois
  representations}, volume 176 of {\em Annals of Mathematics Studies}.
\newblock Princeton University Press, Princeton, NJ, 2011.
\newblock How one can compute in polynomial time the value of Ramanujan's tau
  at a prime.

\bibitem{BG}
B.~H. Gross.
\newblock A tameness criterion for {G}alois representations associated to
  modular forms (mod {$p$}).
\newblock {\em Duke Math. J.}, 61(2):445--517, 1990.

\bibitem{NK}
N.~M. Katz.
\newblock {$p$}-adic properties of modular schemes and modular forms.
\newblock In {\em Modular functions of one variable, {III} ({P}roc. {I}nternat.
  {S}ummer {S}chool, {U}niv. {A}ntwerp, {A}ntwerp, 1972)}, pages 69--190.
  Lecture Notes in Mathematics, Vol. 350, 1973.

\bibitem{NK1}
N.~M. Katz.
\newblock A result on modular forms in characteristic {$p$}.
\newblock In {\em Modular functions of one variable, {V} ({P}roc. {S}econd
  {I}nternat. {C}onf., {U}niv. {B}onn, {B}onn, 1976)}, pages 53--61. Lecture
  Notes in Math., Vol. 601, 1977.

\bibitem{KW}
C.~Khare and J.-P. Wintenberger.
\newblock Serre's modularity conjecture. {I}.
\newblock {\em Invent. Math.}, 178(3):485--504, 2009.

\bibitem{KW1}
C.~Khare and J.-P. Wintenberger.
\newblock Serre's modularity conjecture. {II}.
\newblock {\em Invent. Math.}, 178(3):505--586, 2009.

\bibitem{KM}
M.~Kisin.
\newblock Modularity of 2-adic {B}arsotti-{T}ate representations.
\newblock {\em Invent. Math.}, 178(3):587--634, 2009.

\bibitem{Li}
W.~C.~W. Li.
\newblock Newforms and functional equations.
\newblock {\em Math. Ann.}, 212:285--315, 1975.

\bibitem{OR}
K.~Ono and N.~Ramsey.
\newblock A mod {$\ell$} {A}tkin-{L}ehner theorem and applications.
\newblock {\em Arch. Math. (Basel)}, 98(1):25--36, 2012.

\bibitem{KR}
K.~A. Ribet.
\newblock Report on mod {$l$} representations of {${\rm Gal}(\overline{\bf
  Q}/{\bf Q})$}.
\newblock In {\em Motives ({S}eattle, {WA}, 1991)}, volume~55 of {\em Proc.
  Sympos. Pure Math.}, pages 639--676. Amer. Math. Soc., Providence, RI, 1994.

\bibitem{JS}
J.-P. Serre.
\newblock Sur les repr\'{e}sentations modulaires de degr\'{e} {$2$} de {${\rm
  Gal}(\overline{\bf Q}/{\bf Q})$}.
\newblock {\em Duke Math. J.}, 54(1):179--230, 1987.

\bibitem{WS}
W.~Stein.
\newblock {\em Modular forms, a computational approach}, volume~79 of {\em
  Graduate Studies in Mathematics}.
\newblock American Mathematical Society, Providence, RI, 2007.
\newblock With an appendix by Paul E. Gunnells.

\end{thebibliography}
\end{document}